\input amstex
\documentstyle{amsppt}
\magnification=\magstep1
\font\chapf=cmbx10 scaled \magstep2

\NoBlackBoxes


\nopagenumbers

\redefine\sp{\operatorname{\,sp}}
\define\sa{\operatorname{sa}}
\define\conv{\operatorname{conv}}
\define\Tr{\operatorname{\,Tr}}

\redefine\i{\operatorname{\,i}}

\define\ux{\underline{x}}
\define\uy{\underline{y}}
\define\uz{\underline{z}}
\define\uu{\underline{u}}
\define\us{\underline{s}}

\define\uh{\underline{h}}

\define\ul{\underline{\lambda}}
\define\uI{\underline{I}}

\overfullrule=0pt
\TagsOnRight

\topmatter  
   
\title{\chapf Jensen's trace  \nolinebreak     inequality\\
{}\\
 in several variables}\endtitle
\rightheadtext{Jensen's Trace Inequality}
\author Frank Hansen \& Gert K. Pedersen \endauthor

\date{$4^{th}$ March, 2003}\enddate

\address {Institute of Economics, University of Cop\-enhagen,
Studiestr\ae{}de 6, DK-1455 Copenhagen K, Denmark\;
\&\; Department of Mathematics,  University of Copenhagen, 
Universitetsparken 5, DK-2100 Copenhagen \O, Denmark} 
\endaddress

\email {frank.hansen \@ econ.ku.dk \;\&\; gkped \@ math.ku.dk}
\endemail

\abstract {For a convex, real function $f$ we 
present a simple proof of the formula 
\smallskip
\centerline {$\Tr (f (\sum_{k=1}^m  a_k^* x_ka_k))  
\le \Tr (\sum_{k=1}^m a_k^*f(x_k)a_k)$,}
\smallskip
\noindent valid for each tuple $(x_1, \dots , x_m)$ of 
symmetric matrices in $\Bbb M_n$ and every unital 
column $(a_1, \dots , a_m)$ of matrices, i\.e\.
$\sum_{k=1}^m a^*_ka_k = \bold 1$. This is the 
standard Jensen trace inequality. If $f\ge 0$ it 
holds also for the unbounded trace on 
$\Bbb B(\frak H)$, where $\frak H$ is an 
infinite-dimensional Hilbert space. We then 
investigate the more general case where $\tau$ 
is a  densely defined, lower semi-continuous 
trace on a $C^*-$algebra $\Cal A$ and $f$ is a 
convex, continuous function of $n$ variables, 
and show that we have the inequality 
\smallskip
\centerline{$\tau\left(f(\sum_{k=1}^m a_k^*\ux_ka_k)\right) \le 
\tau\left(\sum_{k=1}^m a_k^* f(\ux_k)a_k\right)$}
\smallskip
\noindent for every family of {\it abelian} $n-$tuples 
$\ux_k= (x_{1k}, \dots, x_{nk})$, i\.e\. tuples of
self-adjoint elements in $\Cal A$ such that 
$[x_{ik},\, x_{jk}]=0$ for all $i, j$ and $k$, 
where $1\le k\le m$, and every unital $m-$column 
$(a_1, \dots, a_m)$ in $M(\Cal A)$, provided that 
the elements $y_i= \sum_{k=1}^m a_k^* x_{ik}a_k$ 
also form an abelian $n-$tuple. We even establish 
this result for weak* measurable, self-adjoint, 
abelian fields $(x_{it})_{t\in T}, \, 1\le i\le n$,
i\.e\. $[x_{it}, x_{jt}]=0$ for all $i, j$ and 
$t$, and a weak* measurable, unital column field 
$(a_t)_{t\in T}$ in $M(\Cal A)$ paired with any 
trace or trace-like functional $\varphi$, i\.e\. 
one that contains the $n-$tuple (presumed abelian) 
with elements $y_i = \int_T a_t^* x_{it}a_t\, d\mu (t)$ 
in its centralizer. This takes the form of the 
inequality
\smallskip
\centerline{$\varphi\left(f(\int_T a_t^*\ux_ta_t\,d\mu(t))\right)\le
\varphi\left(\int_T a_t^*f(\ux_t)a_t\,d\mu(t)\right)$.} 
\smallskip
\noindent We also study functions of $n$ variables that 
are monotone increasing in each variable, and 
show in two important cases that $\varphi (f(\ux)) 
\le \varphi (f(\uy))$ whenever $\ux = 
(x_1,\dots, x_n)$ and $\uy = (y_1,\dots, y_n)$ 
are abelian $n-$tuples with $x_i \le y_i$ for each
$i$ and $\varphi$ is a trace or a trace-like 
functional.} \endabstract 



\endtopmatter
 
\document

\footnote""{\copyright 2001 by the authors. This paper 
may be reproduced, in its entirety, for non-commercial 
purposes.}

\subhead{1. Introduction}\endsubhead Several 
important concepts in operator theory, in quantum 
statistical mechanics (the entropy, the relative 
entropy and Gibbs' free energy), in electrical 
engineering and in mathematical economics involve 
the trace of a function of a self-adjoint operator. 
This has motivated a considerable amount of 
abstract research about such functions in the last
fifty years. An important subset of questions 
concern the convexity of trace functions with 
respect to their argument, and the generalizations 
of this known as Jensen trace inequalities.

The convexity of the function $x@>>> \Tr (f(x))$, 
when $f$ is a convex function of one variable and 
$x$ is a self-adjoint operator, was known to von 
Neumann, cf\. \cite{{\bf 21}, V.3. p. 390}. An 
early proof for $f(x) = \exp(x)$ can be found 
in \cite{{\bf 29}, 2.5.2}. A proof found by 
E.H. Lieb in the early seventies describes the 
number $\Tr(f(x))$, where $f$ is convex, as a  
supremum (taken over all possible choices of 
orthonormal bases of the Hilbert space) of the 
sum of the values of $f$ at the diagonal elements 
of the matrix for $x$. Obviously, then, this is 
a convex function of $x$. The proof was 
communicated to B. Simon, who used the method  
to give an alternative proof of the second 
Berezin-Lieb inequality in \cite{{\bf 30}, 
Theorem 2.4}, see also \cite{{\bf 31}, 
Lemma II.10.4}. Simon only considers the 
exponential function, but the argument is valid 
for any convex function, cf\. \cite{{\bf 17}, 
Proposition 3.1}. The general case for an 
arbitrary normal trace on a von Neumann algebra 
was established by D. Petz in \cite{{\bf 28}, 
Theorem 4}, using the theory of spectral dominance 
(spectral scale). 

When a convex combination $\sum_{k=1}^m \lambda_kx_k$ 
of matrices (or operators) with coefficients 
$(\lambda_1, \dots ,$ $\lambda_m)$ is replaced by the 
non-commutative version $\sum_{k=1}^m a^*_kx_ka_k$, 
where $(a_1, \dots, a_m)$ is a unital $m-$column, 
i\.e\. an $m-$tuple of matrices (or operators) such 
that $\sum_{k=1}^m a_k^*a_k = \bold 1$, we obtain a 
generalization known as Jensen's operator inequality. 
For an {\it operator convex} function, i\.e\. a 
function $f$ such that $f(\lambda x+(1-\lambda)y)
\le \lambda f(x)+(1-\lambda) f(y)$ for any pair of 
self-adjoint matrices $x$ and $y$ (of arbitrary high 
order), this result was found by the first author in 
\cite {{\bf 6}}, and used by the two of us in 
\cite{{\bf 9}} to give a concise review of L\"o{}wner's 
and Bendat-Sherman's theory of operator monotone and 
operator convex functions. With hindsight we must admit 
that we unfortunately chose the contractive form  
$f(a^*xa)\le a^*f(x)a$ for $a^*a\le \bold 1$, this being 
the seemingly most attractive version at the time. 
However, this necessitated the further conditions that 
$0\in I$ and $f(0)\le 0$, conditions that have haunted 
the theory since then, and which become a real obstacle 
when we pass to several variables. The Jensen inequality 
for a trace on a  von Neumann algebra and an arbitrary 
convex function $f$ was found by Brown and Kosaki in 
\cite{{\bf 5}}, still in the contractive version. 
Elementary proofs of these results can now be found in 
\cite{{\bf 11}}.

We begin the paper with the simple proof of the full 
Jensen trace inequality for matrices taken from \cite{{\bf 11}}, 
which uses ideas from Lieb's proof mentioned above. Although 
this result follows from the more general theorem later on in 
the paper we feel that an elementary proof of the most applicable 
version would be a convenience for the (not too specialized) 
reader. Also, the simple proof contains all the basic ideas in 
the more elaborate versions and thus makes it easier to grasp 
these.

\bigskip

\proclaim{2. Theorem} If $f\colon I @>>>\Bbb R$ is a convex 
function defined on an interval $I$ the inequality 
$$
\Tr \left(f\left(\sum_{k=1}^m a_k^*x_ka_k\right)\right)\le
\Tr \left(\sum_{k=1}^m a_k^*f(x_k)a_k \right)
\tag{1}
$$
holds for each $m$-tuple of self-adjoint $n\times n$ 
matrices $(x_1,\dots , x_m)$ with spectra in $I$, every 
unital $m$-tuple $(a_1,\dots,a_m)$ of $n\times n$ matrices
and all natural numbers $n$.
\endproclaim

\demo{Proof} Let $x_k = \sum_{\sp (x_k)} \lambda E_k(\lambda)$ 
denote the spectral resolution of $x_k$ for 
$1\le k \le m$. Thus, $E_k(\lambda)$ is the spectral 
projection of $x_k$ on the eigenspace corresponding 
to $\lambda$ if $\lambda$ is an eigenvalue for $x_k$\,; 
otherwise $E_k(\lambda)=0$. For each unit vector $\xi$ 
in $\Bbb C^n$ define the (atomic) probability measure 
$$
\mu_\xi(S)=\left(\sum_{k=1}^m a_k^* E_k(S)a_k\xi\,\bigg\vert\,\xi\right) = 
\sum_{k=1}^m \left(E_k(S)a_k\xi\mid a_k\xi\right)
\tag{2}
$$
for any (Borel) set $S$ in $\Bbb R$. Note now that
if $y=\sum_{k=1}^m a_k^* x_k a_k$ then 
$$
\aligned
&(y\xi|\xi)=\left(\sum_{k=1}^m 
a_k^* x_k a_k\xi\,\bigg\vert\,\xi\right)\\
=\;&\left(\sum_{k=1}^m\sum_{\sp (x_k)} 
\lambda E_k(\lambda) a_k\xi\,\bigg\vert\, a_k\xi\right)
=\int\lambda\, d\mu_\xi(\lambda).
\endaligned
\tag{3}
$$
If a unit vector $\xi$ is an eigenvector for $y$ then the
corresponding eigenvalue is $(y\xi|\xi)$ and $\xi$ is also 
an eigenvector for $f(y)$ with correponding eigenvalue 
$(f(y)\xi|\xi)=f((y\xi|\xi))$. In this case we therefore have 
$$
\aligned
&\left(f\left(\sum_{k=1}^m a_k^* x_k
a_k\right)\xi\,\bigg\vert\,\xi\right)=(f(y)\xi|\xi)
=f((y\xi|\xi))\\
=\;&f\left(\int\lambda\, d\mu_\xi(\lambda)\right)
\le \int f(\lambda)\, d\mu_\xi(\lambda)\\
=\;&\sum_{k=1}^m\left(\sum_{\sp (x_k)} 
f(\lambda)E_k(\lambda) a_k\xi
\,\bigg\vert\, a_k\xi\right) 
=\sum_{k=1}^m\left(a_k^* f(x_k) a_k\xi\mid\xi\right),
\endaligned
\tag{4}
$$
where we used (3) and the convexity of $f$ -- in form of the 
usual Jensen inequality -- to get the inequality in (4).

The result in (1) now follows by summing over an 
orthonormal basis of eigenvectors for $y$.
\hfill$\square$
\enddemo

\bigskip

\subhead{3. Spectral Theory in Several Variables}\endsubhead 
The really new problems start when we consider a 
function $f(\ul)$ of $n$ real variables (with $\ul = 
(\lambda_1, \lambda_2, \dots , \lambda_n)$). Naturally, 
we wish to replace the real variables $\lambda_j$  by 
self-adjoint operators $x_j$  as in the one-variable 
case. An immediate problem that now arises is how to 
define $f(\ux)$ in this case. The spectral theorem which 
was used in the one-variable case fails here unless the 
$x_j$'s commute with one another. This means that the 
largest domain of definition for $f$ is the set of 
{\it abelian} $n-$tuples in $\Bbb B(\frak H)$, 
i\.e\. tuples $\ux = (x_1, \dots, x_n)$ of self-adjoint 
elements such that $[x_i,\,x_j]=0$ for all $i$ and $j$. 

For functions of two variables the spectral theory of 
abelian tuples (pairs) is equal to spectral theory 
for normal, instead of self-adjoint operators. [As long 
as we consider only continous and not differentiable 
functions, a complex function is just a function of two 
real variables !] This theory is markedly more difficult 
than the one variable case, in particular because the set 
of normal operators has no linear structure.

To be more specific,  consider a $C^*-$algebra 
$\Cal A$ of operators on some Hilbert space 
$\frak H$. For each interval $I$ let $\Cal A_{\sa}^{I}$ 
denote the convex set of self-adjoint elements in 
$\Cal A$ with spectra contained in $I$. If 
$\uI=I_1\times\cdots\times I_n\subset\Bbb R^n$ and $f$ 
is a continuous function on $\uI$ we can for each 
abelian $n-$tuple $\ux=\{x_1,\cdots,x_n\}$ in 
$\bigoplus \Cal A_{\sa}^{I_i}$ define an element $f(\ux)$ 
in $\Cal A$. To see this, let $x_i=\int \lambda dE_i(\lambda)$ 
be the spectral resolution of $x_i$ for $1\le i \le n$. 
Since the $x_i$'s commute, so do their spectral measures. 
We can therefore define the product spectral measure $E$ 
on $\uI$ by $E(S_1\times\cdots\times S_n)=
E_1(S_1)\cdots E_n(S_n)$, and then write   
$$ 
f(\ux)=\int f(\ul)\, dE(\ul) = 
\int f(\lambda_1,\cdots,\lambda_n)\,
dE (\lambda_1,\cdots,\lambda_n).
\tag{5}
$$
Of course, if $f$ is a polynomial in the variables 
$\lambda_1, \dots ,\lambda_n$ we simply find $f(\ux)$ 
by replacing each $\lambda_i$ with $x_i$. The map  
$f\to f(\ux)$ so obtained is a $^*-$homomorphism of 
$C(\uI)$ into $\Cal A$ and generalizes the ordinary
spectral mapping theory for a single (self-adjoint) operator. 
In particular, the support of the map (the smallest closed 
set $S$ such that $f(\ux)=0$ for every function $f$ that 
vanishes on $S$) may be regarded as the ``joint spectrum'' 
of the elements $x_1,\dots , x_n$. In Gelfand language the
commutative unital $C^*-$subalgebra generated by the 
$x_i$'s is $^*-$isomorphic to $C(S)$.

\subhead{4. Convexity in Several Variables}\endsubhead 
The set of abelian $n-$tuples in $\Bbb B(\frak H)$ is 
obviously not a convex set, so at first glance it makes little 
sense to discuss convexity properties of the operator 
function $\ux @>>> f(\ux)$.  We shall therefore consider 
abelian tuples $\ux$ and $\uy$ that are 
{\it compatible}, which by definition means that the line 
segment between them also consists of abelian tuples. 
It is easily seen that this happens precisely when 
$$
[x_i,\,y_j] = [x_j,\,y_i]\quad \text{for all}\; i\; \text{and}\; j.
\tag{6}
$$ 
Now we can meaningfully ask whether 
$\Tr (f(\lambda\ux+(1-\lambda)\uy)\le\Tr(\lambda f(\ux)+
(1-\lambda)f(\uy))$ when $f$ is a convex function. 

Note from (6) that if $\{\ux_1,\dots, \ux_m\}$ is a set of 
pairwise compatible, abelian $n-$tuples, then any linear 
combination $\sum_{k=1}^m \lambda_k\ux_k$ is again an abelian 
$n-$tuple compatible with all the $\ux_i$'s, so that the set 
$\conv \{ \ux_1, \ux_2 \dots , \ux_m\}$ is a convex domain 
for the operator function $f$. This also means that any set 
$\Cal S_o$ of pairwise compatible,  abelian $n-$tuples 
in a $C^*-$subalgebra $\Cal A$ of $\Bbb B(\frak H)$ is 
contained in a maximal set $\Cal S$, which by necessity must 
be a closed, linear subspace of $\Cal A_{\sa}^n$. One may 
wonder how such maximal sets look like, and a few experiments 
show that the variety is wide. Let $\Cal C$ be a commutative 
$C^*-$subalgebra of $\Cal A$ such that $\Cal C'' = \Cal C$ 
(where $'$ denotes relative commutant). For example, $\Cal C$ 
could be the the center of $\Cal A$ (in which case
$\Cal C'=\Cal A$), or it could be any maximal abelian 
$C^*-$subalgebra of $\Cal A$ (in which case $\Cal C' =\Cal C$).
Note though, that the condition $\Cal C'' = \Cal C$ means that 
$\Cal C$ always contains the center of $\Cal A$. Now fix a 
non-zero vector $(\varepsilon_1, \dots, \varepsilon_n)$ in 
$\Bbb R^n$ and define
$$
\Cal S  =\{\ux = (\varepsilon_1 x+c_1, \dots , \varepsilon_n x+c_n) 
\mid x\in \Cal C_{\sa}'\,,\; c_i\in \Cal C_{\sa}\}\,.
\tag{7}
$$
Then it is easy to check that $\Cal S$ is a maximal set of 
pairwise compatible, abelian $n-$tuples in $\Cal A$.

The more useful examples occur, however, at the other 
extreme of the situation above. We assume that the 
$C^*-$algebra $\Cal A$ comes equipped with a set of 
pairwise commuting $C^*-$subalgebras $\Cal A_1, \dots ,
\Cal A_n$. Then the subspace 
$$
\bigoplus_{i=1}^n (\Cal A_i)_{\sa} = \{\ux =(x_1, \dots, x_n)
\mid x_i\in (\Cal A_i)_{\sa}\}
\tag{8}
$$
consists of pairwise compatible, abelian $n-$tuples; and 
under the mild extra condition that each $\Cal A_i$ equals 
the relative commutant in $\Cal A$ of the $C^*-$algebra 
generated by the $\Cal A_j$'s for $j\ne i$, (i\.e\. 
$(\bigcup_{j\ne i} \Cal A_j)' =\Cal A_i$) the space is 
also maximal. This condition may be achieved by replacing 
in turn each of the algebras $\Cal A_i$ by 
$(\bigcup_{j\ne i} \Cal A_j)'$. 

This frame applies readily to the seminal situation 
where $\Cal A=\Cal A_1\otimes\cdots\otimes \Cal A_n$ in 
$\Bbb B (\frak H)$. Indeed, most authors that have 
considered operator functions of several variables have
followed Kor\'a{}nyi's lead and used the functions only 
on tensor products, cf\. \cite{{\bf 14}}.    

In the setting of compatible, abelian tuples we are 
going to replace the trace $\Tr$ on the Hilbert space 
by a densely defined, lower semi-continuous trace  $\tau$ 
on an abstract $C^*-$algebra $\Cal A$; i\.e\. a functional 
defined on the set $\Cal A_+$ of positive elements with 
values in $[0, \infty]$, such that $\tau (x^*x)  =\tau (xx^*)$ 
for all $x$ in $\Cal A$. Thus we shall consider the function
$\ux @>>> \tau (f(\ux))$ on a set of compatible,  abelian 
$n-$tuples in $\Cal A_{\sa}$. 

Some of our results have appeared in more primitive versions 
before. The tracial convexity of the function $\ux @>>> f(\ux)$ 
on the space of $n-$tuples in $\bigoplus_{i=1}^n (\Cal A_i)_{\sa}$ 
with values in $\bigotimes_{i=1}^n \Cal A_i$ was proved by the 
first author for matrix algebras in \cite{{\bf 8}}. His result 
was extended to general operator algebras and traces by the 
second author in \cite{{\bf 27}}. Both proofs rely on Fr\'e{}chet 
differentiability and somewhat intricate manipulations with 
first and second order differentials. It was then realized 
by Lieb that his proof, mentioned above, could be extended 
to the case of several variables with only marginal changes, 
and the improved version appeared in \cite{{\bf 18}}. The 
present version generalizes and subsumes the previous papers. 
In particular we show that the function $\ux @>>>\tau(f(\ux))$  
is convex on any set of the form $\conv \{\ux_1, \dots, \ux_m\}$,   
where the $\ux_k$'s are pairwise compatible,  abelian 
$n-$tuples in $\Cal A_{\sa}$.

\bigskip

\subhead{5. Measurable Fields of Operators}\endsubhead 
Let $\Cal A$ be a (separable) $C^*-$algebra of operators on some 
(separable) Hilbert space $\frak H$ and $T$ a locally compact 
metric space equipped with a Radon measure $\mu$. We say that a field
$(a_t)_{t\in T}$ of operators in the multiplier algebra $M(\Cal A)$ of
$\Cal A$,  i\.e\. the $C^*-$algebra of elements $a$ in $\Bbb B(\frak
H)$ such that  $x\Cal A +\Cal A x \subset \Cal A$,  is {\it weak*
  measurable} if each function  $t @>>> \varphi (a_t)$,  where
$\varphi \in \Cal A^*$, is $\mu-$measurable.  It is worth noticing
that $(a_t)_{t\in T}$ is weak*  measurable if (and only if) for each
vector $\xi$ in $\frak H$  the function $t\to a_t\xi$ is weakly
(equivalently strongly) measurable (because the set  of linear
combinations of vector functionals is weak* dense in  $\Cal A^*$). It
follows that if both $(a_t)_{t\in T}$ and  $(b_t)_{t\in T}$ are weak*
measurable fields then also  $(a_tb_t)_{t\in T}$ is a measurable field.  

If the function $t @>>> \varphi (a_t)$ is integrable for 
all states $\varphi$ and $\int_T |\varphi(a_t)|\,d\mu(t)\le \gamma$ 
for some constant $\gamma$, in particular if the function 
$t@>>> \Vert a_t\Vert$ is integrable, there is a 
unique element in $M(\Cal A)$, designated by $\int_T a_t\,d\mu(t)$, 
such that  
$$
\varphi \left(\int_T a_t\, d\mu (t)\right) = 
\int_T \varphi (a_t)\, d\mu (t)  \qquad \varphi \in \Cal A^*,
\tag{10}
$$
cf\. \cite{{\bf 26}, 2.5.15}. We say in this case that 
the field $(a_t)_{t\in T}$ is {\it integrable}. If all the 
$a_t$'s belong to $\Cal A$  then also  $\int_T a_t\,d\mu(t)$ belongs 
to $\Cal A$. If the weak* measurable field $(a^*_ta_t)_{t\in T}$ 
is integrable with integral $\bold 1$ we say that 
$(a_t)_{t\in T}$ is a {\it unital column field}. 

\bigskip

\subhead{6. Final Notations}\endsubhead Consider now 
an $n-$tuple of weak* measurable, bounded fields 
$(x_{it})_{t\in T}$, each consisting of self-adjoint elements in 
$\Cal A$ with spectra in some fixed interval $I_i$, and assume 
that $[x_{it},\, x_{jt}]=0$ for all $i, j$ and $t$. Thus each vector 
$\ux_t =(x_{1t}, \dots, x_{nt})$ is an abelian $n-$tuple. 
Furthermore, consider a unital column field 
$(a_t)_{t\in T}$ in $M(\Cal A)$, i\.e\. 
$\int_T a_t^*a_t\, d\mu (t) = \bold 1$. Assume finally 
that the elements $y_i = \int_T a_t^*x_{it}a_t\, d\mu (t)$ 
in $\Cal A$ form an abelian $n-$tuple.

The commutation condition above for the $y_i$'s depends on 
intricate relations  between the two measurable fields 
$(x_{it})_{t\in T}$ and $(a_t)_{t\in T}$. It is, 
however, satisfied if the fields satisfy the following 
extension of the commutativity condition in (5):
$$
[a^*_tx_{it}a_t,\, a^*_sx_{js}a_s] = 
[a^*_tx_{jt}a_t,\, a^*_sx_{is}a_s]\quad \text{for all} 
\; s\;\text{and}\;t.
\tag{11}
$$ 
Thus in particular if $x_{it}a_ta^*_sx_{js}=x_{jt}a_ta^*_sx_{is}$ 
for all $s$ and $t$.
 
For ease of notation we shall write 
$\uI = I_1 \times \cdots \times I_n$ and 
$\ul = (\lambda_1, \dots ,\lambda_n)$ if 
$\ul \in \uI$. Moreover, we regard the vector 
space of $n-$tuples in $\Cal A$ as a bimodule over 
$M(\Cal A)$ and write $\ux_t=(x_{1t},\dots, x_{nt})$ 
and $a^*_t\ux_ta_t=(a^*_tx_{1t}a_t,\dots, a^*_tx_{nt}a_t)$, 
so that $\uy=(y_1,\dots, y_n)=\int_T a^*_t\ux_ta_t\, d\mu (t)$. 

We finally recall that the {\it centralizer} of a positive 
functional $\varphi$ on $\Cal A$ is the $C^*-$subalgebra 
$\Cal A^{\varphi} =  \{y\in \Cal A\mid \forall x \in \Cal A
\, \colon \, \varphi(xy)=\varphi(yx)\}$. If $\varphi$ is unbounded, 
but lower semi-continuous on $\Cal A_+$ and finite on 
the minimal dense ideal $K(\Cal A)$ of $\Cal A$, we define 
$\Cal A^{\varphi}=\{y\in \Cal A\mid\forall x\in K(\Cal A)
\,\colon\,\varphi(xy)=\varphi(yx)\}$.

\bigskip

\proclaim{7. Theorem} Let $(\ux_t)_{t\in T}$ be a bounded, 
weak* measurable field of abelian $n-$tuples in a 
$C^*-$algebra $\Cal A$, with $\sp (x_{it}) \subset I_i$ for 
$1\le i\le n$, and let $(a_t)_{t\in T}$ be a unital column field in 
$M(\Cal A)$ such that the elements $y_i =
\int_T a^*_tx_{it}a_t\, d\mu(t)$ form an abelian $n-$tuple. 
Then for each continuous,  convex function $f$ defined on the 
cube $\uI = I_1\times \cdots \times I_n$ in $\Bbb R^n$ and 
every positive functional $\varphi$ that contains the $y_i$'s 
in its centralizer $\Cal A^{\varphi}$, i\.e\. $\varphi (xy_i) = 
\varphi (y_ix)$ for all $x$ in $\Cal A$ and every $y_i$, we have 
the inequality:
$$
\varphi\left(f\left(\int_T a_t^*\ux_t a_t\,d\mu (t)\right)\right)
\le \varphi\left(\int_T a^*_t f(\ux_t)a_t\,d\mu (t)\right).
\tag{12}
$$

If $\varphi$ is unbounded, but lower semi-continuous on 
$\Cal A_+$ and finite on the minimal dense ideal 
$K(\Cal A)$ of $\Cal A$, the result still holds if $f\ge 0$, 
even though the function may now attain infinite values. 
\endproclaim

\demo{Proof} Let $\Cal C = C_o (S)$ denote the commutative 
$C^*-$subalgebra of $\Cal A$ generated by $y_1, \dots , y_n$, 
and let $\mu_\varphi$ be the finite Radon measure on 
the locally compact, metric space $S$ defined, via the 
Riesz representation theorem, by 
$$
\int_S y(s)\,d\mu_{\varphi} (s)=\varphi(y)\qquad y\in \Cal C=C_o(S).
\tag{13}
$$
Since for all $(x,y)$ in $M(\Cal A)_+\times \Cal C_+$ we have $\varphi(xy)=
\varphi(y^{1/2}xy^{1/2})$ it follows that
$$
0\le \varphi(xy) \le \Vert x \Vert \varphi(y).
\tag{14}
$$
Consequently the functional $y\to\varphi(xy)$ on $\Cal C$ 
defines a Radon measure on $S$ dominated by a multible 
of $\mu_{\varphi}$, hence determined by a unique element 
$\Phi(x)$ in $L^\infty_{\mu_\varphi} (S)$. By linearization 
this defines a conditional expectation 
$\Phi\colon M(\Cal A)\to L^\infty_{\mu_\varphi}(S)$ 
(i\.e\. a positive, unital module map) such that 
$$
\int_S y(s)\Phi(x)(s)\,d\mu_\varphi (s)=\varphi(yx),\qquad y\in \Cal C 
\quad x\in M(\Cal A).
\tag{15}
$$
Inherent in this formulation is the fact that if 
$y\in \Cal C = C_o(S)$,  then $\Phi(y)$ is the natural image of 
$y$ in $L^\infty_{\mu_\varphi}(S)$. In particular, 
$y(s) = \Phi(y)(s)$ for almost all $s$ in $S$. 

Observe now that since the $C^*-$algebra $C_o(\uI)$ is 
separable we can for almost every $s$ in $S$ define a 
Radon measure $\mu_s$ on $\uI$ by
$$
\int_{\uI} g(\ul)\, d\mu_s (\ul) = \Phi\left( \int_T a^*_t
g(\ux_t)a_t\,d\mu(t)\right)(s) \qquad g\in C_o(\uI).
\tag{16}
$$
As $\int_Ta^*_ta_t\,d\mu(t)=\bold 1$  this is actually 
a probability measure. 

If we put $g_i(\ul)= \lambda_i$ for $1\le i \le n$ then 
$$
\int_{\uI} g_i(\ul)\,d\mu_s(\ul) = \Phi\left(\int_T a^*_t x_{it} a_t\,
d\mu(t)\right)(s) = \Phi(y_i) (s) = y_i(s).
\tag{17}
$$
Since $y_i \in \Cal C$ for all $i$ we get by (17) -- using 
the convexity of $f$ in form of the standard Jensen 
inequality -- that
$$
\aligned
f(\uy)(s)\;&=f(\uy(s))=f(y_1(s),\dots, y_n(s))\\
&=f\left(\int_{\uI} g_1(\ul)\,d\mu_s(\ul), \dots, \int_{\uI} g_n(\ul)\,
d\mu_s(\ul)\right) \\
&\le\ \int_{\uI} f\left(g_1(\ul), \dots , g_n(\ul)\right)\, d\mu_s(\ul)\\
&= \int_{\uI} f(\ul)\,d\mu_s(\ul) =   
\Phi\left(\int_T a^*_tf(\ux_t)a_t\,d\mu(t)\right)(s).
\endaligned
\tag{18}
$$
Integrating over $s$ now gives the desired result:
$$
\aligned
&\varphi(f(\uy)) = \int_S f(\uy)(s)\, d\mu_\varphi(s)\\ 
\le\;&\int_S \Phi\left( \int_T  a^*_t f(\ux_t)a_t\, d\mu(t)\right) (s)\, 
d\mu_\varphi(s)\\ 
=\;& \int_T\int_S \Phi \left( a^*_tf(\ux_t)a_t \right)(s)\, 
d\mu_\varphi (s)\,d\mu(t)\\
=\;&\int_T \varphi\left( a^*_tf(\ux_t)a_t\right)\, d\mu(t)
= \varphi\left( \int_T a^*_t f(\ux_t)a_t\, d\mu(t)\right). 
\endaligned 
\tag{19}
$$
 
Having proved the finite case, let us now assume that 
$\varphi$ is unbounded, but lower semi-continuous on 
$\Cal A_+$ and finite on the minimal dense ideal 
$K(\Cal A)$ of $\Cal A$ . Such functionals were 
termed $C^*-integrals$ in \cite{{\bf 23}} and 
\cite{{\bf 24}}. This -- by definition -- means that 
$\varphi(x) < \infty$ if $x\in \Cal A_+$ and $x=xe$ 
for some $e$ in $\Cal A_+$, because $K(\Cal A)$ is 
the hereditary $^*-$subalgebra of $\Cal A$ 
generated by such elements, cf\. \cite{{\bf 25}, 5.6.1}. 
Restricting $\varphi$ to $\Cal C$ we therefore obtain a 
unique Radon measure $\mu_\varphi$ on $S$ such that    
$$
\int_S y(t)\,d\mu_\varphi (t) = \varphi (y) \qquad y\in \Cal C.
\tag{20}
$$
Inspection of the proof above now shows that the 
Jensen trace inequality still holds if only $f \ge 0$, 
even though $\infty$ may now occur in the  inequality. 
\hfill$\square$
\enddemo

\bigskip

\subhead{8. Remarks}\endsubhead The second condition in 
Theorem 7, that the elements $y_i$ are mutually commuting, 
is not easy to verify. There are, however, a few cases 
that can be handled with ease. In the first we simply 
set $n=1$, so that we obtain the one-variable extension 
of Theorem 2. This is done in Corollary 9. In the second 
case we let each $a_t$ be a positive scalar and set 
$a^*_ta_t = \lambda (t)$. Then if $[x_{it},\, x_{js}]
=[x_{jt},\, x_{is}]$ for all $i,j,s$ and $t$ (so, in 
particular $[x_{it},\, x_{jt}]=0$), the elements 
$y_i=\int_T x_{it}\lambda (t)\,d\mu (t)$ will form a 
abelian $n-$tuple. Thus in Corollary 11 we obtain an 
extremely strong version of the convexity of the trace 
function, proved in weaker forms in \cite{{\bf 8, 27, 18}}. 

\bigskip

\proclaim{9. Corollary} For each convex, continuous 
function $f$ on an interval $I$, every bounded, weak* 
measurable field $(x_t)_{t\in T}$ in
$\Cal A_{\sa}^I$ and every unital column field 
$(a_t)_{t \in T}$ in $M(\Cal A)$ we have the inequality
$$
\varphi\left(f\left(\int_T a^*_tx_ta_t\, d\mu (t)\right)\right) \le
\varphi\left(\int_T a^*_t f(x_t)a_t\,d\mu(t)\right) 
\tag{21}
$$
for every positive functional $\varphi$ that contains the 
element $y=\int_Ta^*_tx_ta_t\, d\mu(t)$ in its centralizer. 

If $\varphi$ is unbounded, but lower semi-continous and 
finite on the minimal dense ideal $K(\Cal A)$ of $\Cal A$, the result 
still holds if $f\ge 0$. \hfill$\square$
\endproclaim

For continuous fields this result was proved in \cite{{\bf 11},
  Theorem 4.1}. 
\bigskip

\proclaim{10. Corollary} For each convex, continuous function 
$f$ on a cube $\uI=I_1\times\cdots I_n$ in $\Bbb R^n$, each
probability measure $\mu$ on a locally compact,  metric space 
$T$ and every $n-$tuple of bounded, weak* measurable fields 
$(x_{it})_{t\in T}$, where $x_{it} \in \Cal A_{\sa}^{I_i}$ for 
$1\le i \le n$, such that $[x_{it},\, x_{js}]=[x_{jt},\, x_{is}]$ for
all $i, j, s$ and $t$ we have
$$
\varphi\left(f\left(\int_T \ux_t\, d\mu (t)\right)\right) \le
\varphi\left(\int_T (f(\ux_t))\,d\mu (t)\right),
\tag{22}
$$ 
for every positive functional $\varphi$ on $\Cal A$ that contains 
the elements $y_i = \int_T x_{it}\, d\mu (t)$ in its centralizer. 
\hfill$\square$ 
\endproclaim

\bigskip

Specializing to convex combinations (discrete probability measures)
and traces we obtain the following version of Corollary 10:

\proclaim{11. Corollary} For each convex, continuous function 
$f$ on a cube $\uI=I_1\times\cdots I_n$ in $\Bbb R^n$, and every trace
$\tau$ on a $C^*-$algebra $\Cal A$  the function
$$
\ux @>>>\tau(f(\ux))
\tag{23}
$$
is convex on the set of compatible pairs of abelian 
$n-$tuples $\ux = (x_1,\dots, x_n)$ in $\Cal A$ with 
$\sp (x_i)\in I_i$ for all $i$. \hfill$\square$ 
\endproclaim

\bigskip

The condition in Theorem 7 that the elements $x_{it}$ and 
$x_{jt}$ commute mutually is also rather awkward. An easy and 
important solution to this problem is to assume from  
the outset that the $C^*-$algebra $\Cal A$ comes equipped 
with mutually commuting $C^*-$subalgebras 
$\Cal A_1,\dots, \Cal A_n$ and then require that $x_{it}\in \Cal A_i$
for all $i$ and $t$. Now the domain of definition of 
$f$ is the convex set $\bigoplus_{i=1}^n (\Cal A_i)_{\sa}^{I_i}$ 
and we can state the Jensen trace inequality for $f$ 
in ordinary terms.

\proclaim{12. Corollary} For each convex, continuous function 
$f$ on a cube $\uI=I_1\times\cdots\times I_n$ in $\Bbb R^n$
and every $C^*-$algebra $\Cal A$ with mutually commuting 
$C^*-$subalgebras $\Cal A_1, \dots \Cal A_n$ we have the 
inequality 
$$
\varphi\left(f\left(\int_T a^*_t\ux_ta_t\,d\mu(t)\right)\right)
\le \varphi \left(\int_T a^*_tf(\ux_t) a_t\, d\mu (t)\right)
\tag{24}
$$
for each bounded, weak* measurable field $(\ux_t)_{t\in T} = 
\left((x_{1t})_{t\in T},\dots, (x_{nt})_{t\in T}\right)$  in  
$\bigoplus_{i=1}^n (\Cal A_i)^{I_i}_{\sa}$, and every unital 
column field $(a_t)_{t\in T}$ in $M(\Cal A)$, provided 
that the elements $y_i=\int_T a^*_tx_{it}a_t\,d\mu(t)$ 
form an abelian $n-$tuple and the functional $\varphi$ 
contains these elements in its centralizer. 

If $\varphi$ is unbounded, but lower semi-continous and finite 
on the minimal dense ideal $K(\Cal A)$ of $\Cal A$, the result 
still holds if $f\ge 0$. \hfill$\square$
\endproclaim

This result generalizes both \cite{{\bf 11}, Theorem 4.1} and
\cite{{\bf 27}, Theorem 2}.

\bigskip

In the next case let the parameter space be 
$\Bbb N_n\times T$, where $\Bbb N_n$ denotes the finite 
subset $\{1,2, \dots, n\}$. So instead of the index $t$ 
we now have the double index $(j, t)$. We then assume that
$x_{ijt} = x_{it}$ if $j=i$ and that $x_{ijt}=0$ if $j\ne i$. 
Furthermore we assume that $a_{jt}\in M(\Cal A_j)$ for all $j$, 
so that we have the elements $b_j =\int_T a^*_{jt}a_{jt}\, 
d\mu (t)$ in $M(\Cal A_j)$ with  $\sum_{j=1}^n b_j =\bold 1$. 
Note now that 
$$
y_i = \sum_{j=1}^n \int_T a^*_{jt}x_{ijt}a_{jt}\, d\mu (t) = 
\int_T a^*_{it}x_{it}a_{it}\,d\mu(t) \in \Cal A_i,
\tag{25}
$$
so the commutativity condition is trivially satisfied. Consequently we
have the following result:

\proclaim{13. Corollary} For each convex, continuous function 
$f$ on a cube $\uI=I_1\times \cdots \times I_n$, where $0\in I_i$ 
for each $i$, every $n-$tuple of bounded, weak* measurable 
fields $(x_{it})_{t\in T} \subset (\Cal A_i)_{\sa}^{I_i}$ and 
every $n-$tuple of integrable column fields $(a_{it})_{t \in T}$ in 
$M(\Cal A_i)$ with  $\sum_{i=1}^n\int_Ta_{it}^*a_{it}\,d\mu(t)=\bold 1$
we have the inequality:   
$$
\aligned
&\varphi\left( f\left( \int_T a^*_{1t}x_{1t}a_{1t}\,d\mu(t), 
\dots , \int_T a^*_{nt}x_{nt}a_{nt}\,d\mu(t) \right)\right) \\
\le\;&\varphi \left(\sum_{i=1}^n \int_T a^*_{it} f(0, \dots , x_{it},
\dots , 0)a_{it}\, d\mu(t)\right) 
\endaligned
\tag{26}
$$
for every positive functional $\varphi$ that contains the 
elements $y_i=\int_T a^*_{it}x_{it}a_{it}\, d\mu(t)$ in 
its centralizer. \hfill$\square$\endproclaim

\bigskip

In the last case we again use the parameter space 
$\Bbb N_n \times T$ and take $a_{jt}$ in $M(\Cal A_j)$ for 
all $j$, but we now put $x_{ijt} = x_i$ constantly for 
some fixed $x_i$ in $(\Cal A_i)_{\sa}$. Then 
$$
y_i = \sum_{j=1}^n \int_T a^*_{jt}x_ia_{jt}\, d\mu (t) = 
\int_T a^*_{it}x_ia_{it}\,d\mu(t) + (\bold 1 -b_i)x_i \in \Cal A_i,
\tag{27}
$$
so again we have the desired relations. This gives the following
result: 

\proclaim{14. Corollary} For each convex, continuous 
function $f$ on a cube $\uI$, every $n-$tuple $\ux$ with elements 
$x_i$ in $ (\Cal A_i)_{\sa}^{I_i}$ and every $n-$tuple of integrable 
column fields $(a_{it})_{t\in T}$ in $M(\Cal A_i)$,  with 
$\sum_{i=1}^n b_i = \bold 1$, where 
$b_i = \int_T a_{it}^*a_{it} \,d\mu(t)$, we have the inequality:  
$$
\align
&\varphi\left(f\left(\int_T a^*_{1t}x_1a_{1t}\,d\mu(t)+
 (\bold 1-b_1)x_1,\dots, \int_T a^*_{nt}x_na_{nt}\,d\mu(t)+
 (\bold 1 -b_n)x_n\right)\right) \\
\le\; &\varphi\left(\sum_{i=1}^n\int_T
a^*_{it}f(x_1,\dots,x_n)a_{it}\,d\mu(t)\right)    \tag{28}
\endalign
$$
for every  positive functional $\varphi$ that contains 
the elements $y_i=\int_Ta^*_{it}x_ia_{it}\,d\mu(t)+
(\bold 1-b_i)x_i$ in its centralizer. \hfill$\square$ 
\endproclaim

\bigskip
 
\subhead{15. Monotonicity}\endsubhead We conclude the 
paper with some results about monotonicity of operator 
functions under a trace or a trace-like functional. 
The tendency is that if $f$ is monotone increasing 
in each variable and $\ux=(x_1,\dots,x_n)$ and 
$\uy=(y_1,\dots,y_n)$ are abelian tuples, so that we 
have a chance to define $f(\ux)$ and $f(\uy)$, then 
$\varphi(f(\ux))\le \varphi (f(\uy))$ if only 
$x_i\le y_i$ for all $i$. This result may or may not be 
true in general. We can prove it when $f$ is either convex 
or concave, or when $\ux$ and $\uy$ are compatible.
 
\bigskip

\proclaim{16. Theorem} Let $f\colon \uI @>>> \Bbb R$ 
be a continuous function on a cube 
$\uI=[\alpha_1, \beta_1]\times \cdots \times [\alpha_n, \beta_n]$ 
in $\Bbb R^n$, and assume that $f$ is monotone 
increasing in each variable. Then for any two abelian 
$n-$tuples $(x_1,\dots,x_n)$ and $(y_1,\dots, y_n)$ 
of self-adjoint elements in a $C^*-$algebra $\Cal A$ with 
$\alpha_i\bold 1 \le x_i \le y_i \le \beta_i\bold 1$ 
for all $i$ we have the inequality 
$$
\varphi (f(\ux))\le \varphi (f(\uy))
\tag{29}
$$
for any positive functional $\varphi$ on $\Cal A$ that contains 
the elements $x_1,\dots, x_n$ in its centralizer, provided 
that $f$ is also convex. If instead $f$ is concave the 
result holds if the elements $y_1,\dots, y_n$ belong to 
the centralizer of $\varphi$.
\endproclaim

\demo{Proof} Let $\Cal C=C_o(S)$ denote the commutative 
$C^*-$subalgebra of $\Cal A$ generated by the $x_i$'s. 
As in the proof of Theorem 7 we then obtain a Radon
measure $\mu_{\varphi}$ on $S$ and a conditional 
expectation $\Phi\colon M(\Cal A) @>>>L_{\mu_\varphi}^{\infty}(S)$ 
such that
$$
\int_S z(s)\Phi(y)(s)\,d\mu_\varphi(s)
=\varphi(zy), \qquad z\in \Cal C \quad y\in M(\Cal A),
\tag{30}
$$
where $\Phi(z)(s) =z(s)$ almost everywhere on $S$ 
for each $z$ in $\Cal C$.

Since $\uI$ is separable we can for almost every $s$ 
in $S$ define a probability measure $\mu_s$ on $\uI$ 
by the formula
$$
\int_{\uI}g(\ul)\,d\mu_s(\ul)=\Phi(g(\uy))(s)\qquad g\in C(\uI).
\tag{31}
$$
If we set $g_i(\ul)= \lambda_i$ for each $i$, this 
means that
$$
\int_{\uI} g_i(\ul)\, d\mu_s(\ul)= \Phi(g_i(\uy))(s)= \Phi(y_i)(s).
\tag{32}
$$

Assume now that $f$ - in addition to being monotone 
increasing - is also convex on the cube $\uI$. Then, 
using that $f(\ux)\in \Cal C$ it follows that
$$
\aligned
&\Phi(f(\ux))(s)=f(\ux)(s)=f(x_1(s),\dots, x_n(s))\\
=\; &f\left(\Phi(x_1)(s),\dots, \Phi(x_n)(s)\right)
\le f\left(\Phi(y_1)(s),\dots, \Phi(y_n)(s)\right)\\
=\; & f\left(\int_{\uI} g_1(\ul)\,d\mu_s(\ul),\dots, \int_{\uI}
g_n(\ul)\,d\mu_s(\ul)\right)\\
\le\;& \int_{\uI}f\left(g_1(\ul),\dots, g_n(\ul)\right)\,d\mu_s(\ul)
=\int_{\uI} f(\lambda_1,\dots, \lambda_n)\,d\mu_s(\ul)\\
=\;&\int_{\uI} f(\ul)\,d\mu_s (\ul)=\Phi(f(\uy))(s), 
\endaligned 
\tag{33}
$$
where we used the monotonicity of $f$ to obtain the first 
inequality sign in (33) and the convexity of $f$ -- 
in form of the usual Jensen inequality -- to obtain the 
second inequality sign. Integrating over $s$ now yields the
desired result:
$$
\aligned
&\varphi(f(\ux))=\int_S \Phi (f(\ux))(s)\,d\mu_\varphi(s) \\
\le\;& \int_S \Phi (f(\uy))(s)\,d\mu_\varphi(s) = \varphi (f(\uy)). 
\endaligned
\tag{34}
$$

If on the other hand we assume that $f$ is a concave 
function we simply permute the r\^ oles of the $n-$tuples 
$\ux$ and $\uy$ and let $\Cal C$ denote the $C^*-$subalgebra 
of $\Cal A$ generated by the $y_i$'s. The conditional 
expectation $\Phi\colon M(\Cal A)@>>> L_{\mu_\varphi}^{\infty}(S)$ 
now satisfies that $\Phi (y_i)(s)=y_i(s)$ almost 
everywhere. Similarly we redefine the probability
measures $\mu_s$ by the new formula 
$\int_{\uI} g(\ul)\,d\mu_s(\ul)= \Phi(g(\ux))(s)$, 
so that now
$$
\int_{\uI} g_i(\ul)\, d\mu_s(\ul)= \Phi(g_i(\ux))(s)= \Phi(x_i)(s).
\tag{35}
$$
It follows as in (33) that we have the inequalities
$$
\aligned
&\Phi(f(\ux))(s)= \int_{\uI}f(\ul)\,d\mu_s (\ul)\\
=\;&\int_{\uI} f(\lambda_1,\dots, \lambda_n)\,d\mu_s(\ul)
=\int_{\uI}f\left(g_1(\ul),\dots, g_n(\ul)\right)\,d\mu_s(\ul)\\
\le\;& f\left(\int_{\uI} g_1(\ul)\,d\mu_s(\ul),\dots, \int_{\uI}
g_n(\ul)\,d\mu_s(\ul)\right)\\
=\;& f\left(\Phi(x_1)(s),\dots, \Phi(x_n)(s)\right)
\le f\left(\Phi(y_1)(s),\dots, \Phi(y_n)(s)\right)\\
=\;&f(y_1(s),\dots, y_n(s))=f(\uy)(s)=\Phi(f(\uy))(s),
\endaligned
\tag{36}
$$
where we now used the concavity of $f$ to obtain the 
first inequality in (36). Integrating over $s$ we 
again get the desired inequality (29).
\hfill$\square$
\enddemo

\bigskip

\subhead{17. Remarks}\endsubhead Evidently we may 
combine the two conditions in Theorem 16 to show 
that if $f$ is an increasing function which admits 
a decomposition $f=f_++f_-$, where $f_+$ and $f_-$ 
are both increasing and $f_+$ is convex whereas 
$f_-$ is concave, then $\varphi(f(\ux)) \le 
\varphi(f(\uy))$ if $\ux \le \uy$, provided that  
all the elements $x_1,\dots, x_n$ and 
$y_1,\dots, y_n$ belong to the centralizer of 
$\varphi$. However, such a decomposition, even 
approximately, is not possible in general, not 
even in the one-variable case. The reader may 
check that $\sin (t)$, for $-\pi/2 \le t \le \pi/2$ 
can not be approximated by any function 
$f=f_++f_-$, where $f_+$ is convex and  
$f_-$ is concave, and both are increasing. In fact, 
$\Vert \sin -f\Vert_{\infty} > (2\pi)^{-2}$. 

The simple function $f(s, t)=st$ is neither convex 
nor concave, but increases in each variable on the 
first quadrant. One  proves by direct calculations 
that if $x_1, y_1$ and $x_2, y_2$ are positive
elements in a $C^*-$algebra $\Cal A$ with 
$x_1 \le x_2$, $y_1 \le y_2$ then 
$ \tau (x_1y_1) \le \tau (x_2y_2)$ for every 
trace  $\tau$ on $\Cal A$. The simple argument 
relies on the cyclicity of the trace, which for two 
factors is equivalent to commutativity, but does not 
need the commutator equations $[x_1, y_1]=[x_2, y_2]=0$, 
which we are prepared to insert to get abelian tuples. 
This particular argument fails for three factors, so that we are 
not able to decide whether the function $f(r,s,t)=rst$ 
is an increasing trace function on positive abelian 
triples. 

Despite this setback one may still hope that the 
function $\ux @>>> \tau (f(\ux))$ is increasing 
on the set of abelian $n-$tuples, provided only that 
$f$ is monotone increasing in each variable; at least 
when $\tau$ is a trace or a trace-like functional. 
Our last result, an extension of \cite {{\bf 27}, 
Corollary 5}, shows that this is true  when the two
abelian $n-$tuples are compatible. The proof uses 
the Fr\'e{}chet differential as in \cite{{\bf 10}}. 

\bigskip

\proclaim{18. Proposition} Let $f$ be a continuous 
function on a cube $\uI=[\alpha_1, \beta_1]\times
\cdots \times [\alpha_n, \beta_n]$ in $\Bbb R^n$, 
and assume that $f$ is increasing in each variable. 
Then for any two compatible,  abelian $n-$tuples 
$\ux$ and $\uy$ in a $C^*-$algebra $\Cal A$ that 
satisfy $\alpha_i\le x_i \le y_i\le \beta_i$ for all 
$i$, we have the inequality $\varphi(f(\ux))\le 
\varphi(f(\uy))$ for any positive functional $\varphi$ 
on $\Cal A$ that contains all the elements $x_1,\dots, x_n$ and 
$y_1,\dots, y_n$ in its centralizer.
\endproclaim

\demo{Proof} Put $\uh=\uy-\ux$, and 
let $g(t)=\varphi(f(\ux +t\uh))=
\varphi(f((1-t)\ux+t\uy))$ for $0\le t\le 1$. 
(Note that this is well defined since $\ux$ 
and $\uy$ are compatible.) Then
$$
\varphi(f(\uy))-\varphi(f(\ux))=\int_0^1 g'(t)\,dt,
\tag{37}
$$
provided, of course, that $g$ is differentiable. 
However, working by approximation -- extending $f$ 
to a bounded increasing function on $\Bbb R^n$ and 
convolving it with a suitable approximate unit for 
$L^1(\Bbb R^n)$ like $(\varepsilon/\pi)^{n/2}
\exp (-\varepsilon \us\cdot\us)$ --  we may assume 
that $f$ is extendable to a Schwartz function on 
$\Bbb R^n$, whence $f(\uu)=\int_{\Bbb R^n} 
\exp(\i(\uu\cdot\us))\widehat f(\us)\,d\us$. 
Consequently, with $\uz=\ux+t\uh$, 
$$
\aligned
g'(t)\;&=\lim \, \varepsilon^{-1}
\varphi(f(\uz+\varepsilon\uh)-f(\uz))\\
&=\lim \,\varepsilon^{-1}\int_{\Bbb R^n}
\varphi\left(\exp(\i((\uz+\varepsilon\uh)\cdot\us))-
\exp(\i(\uz\cdot\us))\right)\widehat f(\us)\,d\us.
\endaligned
\tag{38}
$$
By the Dyson expansion of the operator function 
$b @>>> \exp (a+b)$ we have the expression 
$\lim \, \varepsilon^{-1} (\exp (a+\varepsilon b)-\exp (a))= 
\int_0^1 \exp(ra)b\exp((1-r)a)\,dr$, and inserting 
this in (38) we get
$$
\aligned
g'(t)\; &=\int_{\Bbb R^n}\varphi\left(\int_0^1
\exp(\i(\uz\cdot\us)r)\i(\uh\cdot\us)
\exp(\i(\uz\cdot\us)(1-r))
\,dr\right)\widehat f(\us)\,d\us \\
&=\int_{\Bbb R^n} \varphi\left(\exp(\i(\uz\cdot\us))
\i(\uh\cdot\us)\right)\widehat f(\us)\,d\us \\
&=\sum_{k=1}^n\int_{\Bbb R^n} 
\varphi\left(\exp(\i(\uz\cdot\us))h_k\right)\i s_k\widehat f(\us)\,d\us
=\sum_{k=1}^n \varphi(f'_k(\uz)h_k),
\endaligned
\tag{39}
$$
where we used that the element $\uz\cdot\us$, hence also 
$\exp(\i(\uz\cdot\us)(1-r))$, is in the centralizer of $\varphi$.
Since $f'_k\ge 0$ and $h_k\ge 0$ for all $k$ it follows 
that $g'\ge 0$, whence $\varphi(f(\ux))\le \varphi(f(\uy))$ by
(37), as desired.
\hfill$\square$
\enddemo


\enddocument